\documentclass[a4paper,12pt]{article}
  \usepackage[top=2.5cm,bottom=2.5cm,left=2.5cm,right=2.5cm]{geometry}
\usepackage{mathrsfs}

\usepackage{amsmath, epsfig, cite}
\usepackage{amssymb}
\usepackage{amsfonts}
\usepackage{latexsym}
\usepackage{graphicx}

\newtheorem{thm}{Theorem}[section]
\newtheorem{prop}[thm]{Proposition}

\newtheorem{defi}[thm]{Definition}
\newtheorem{lem}[thm]{Lemma}

\newcommand{\qed}{{\hfill\rule{4pt}{7pt}}}
\def\pf{\noindent {\it Proof.} }

\numberwithin{equation}{section}

\makeatletter \@addtoreset{equation}{section} \makeatother

\begin{document}
\rule{0cm}{3.5cm}

\begin{center} {\Large \bf The Maximal Matching Energy of Tricyclic Graphs}
\end{center}
\begin{center}
{{ \small Lin Chen, Yongtang Shi\footnote{The corresponding author.}} \\[3mm]
{\small Center for Combinatorics and LPMC-TJKLC, Nankai University}\\
{\small Tianjin 300071, P.R. China.}\\[3mm]
{\small E-mail: chenlin1120120012@126.com, shi@nankai.edu.cn } \\[2mm]
(Received May 13, 2014)}
\end{center}
\begin{center}
\begin{minipage}{120mm}
\vskip 0.2cm
\begin{center}
{\bf Abstract}
\end{center}
{\small Gutman and Wagner proposed the concept of the matching energy (ME) and pointed
out that the chemical applications of ME go back to the 1970s. Let $G$ be a simple graph
of order $n$ and $\mu_1,\mu_2,\ldots,\mu_n$ be the roots of its matching polynomial. The
matching energy of $G$ is defined to be the sum of the absolute values of $\mu_{i}\
(i=1,2,\ldots,n)$. Gutman and Cvetkoi\'c determined the tricyclic graphs on $n$ vertices
with maximal number of matchings by a computer search for small values of $n$ and by an
induction argument for the rest. Based on this result, in this paper, we characterize
the graphs with the maximal value of matching energy among all tricyclic graphs, and
completely determine the tricyclic graphs with the maximal matching energy. We prove our
result by using Coulson-type integral formula of matching energy, which is similar as
the method to comparing the energies of two quasi-order incomparable graphs.}

\end{minipage}
\end{center}

 \baselineskip=0.30in

\section{Introduction}

In this paper, all graphs under our consideration are finite, connected, undirected
and simple. For more notations and terminologies that will be used in the sequel, we
refer to \cite{graphBondy2008}. Let $G$ be such a graph, and let $n$ and $m$ be the
number of its vertices and edges, respectively. A \emph{matching} in a graph $G$ is
a set of pairwise nonadjacent edges. A matching $M$ is called a $k$-\emph{matching}
if the size of $M$ is $k$. Denote by $m(G,k)$ the number of $k$-matchings of $G$,
where $m(G,1)=m$ and $m(G,k)=0$ for $k>\lfloor\frac{n}{2}\rfloor$ or $k<0$. In
addition, define $m(G,0)=1$. The matching polynomial of graph $G$ is defined as
$$\alpha(G)=\alpha(G,x)=\sum\limits_{k\geq 0}(-1)^km(G,k)x^{n-2k}.$$

In 1977, Gutman \cite{gutman1977} proposed the concept of graph energy. The energy
of $G$ is defined as the sum of the absolute values of its eigenvalues, namely,
\begin{equation*}
E(G)=\sum\limits_{i=1}^{n}|\lambda_i|,
\end{equation*}
where $\lambda_1,\lambda_2,\ldots,\lambda_n$ denote the eigenvalues of $G$. The theory
of graph energy is well developed. The graph energy has been rather widely studied by
theoretical chemists and mathematicians. For details, we refer the book on graph energy
\cite{li&shi&gutman2012} and reviews \cite{gutman2001, gutman2009}. Recently, Gutman and
Wagner \cite{gutman&Wagner2012} defined the matching energy of a graph $G$ based on the
zeros of its matching polynomial \cite{Farrell1979,gutmanmatch1979}.

\begin{defi}
Let $G$ be a simple graph with order $n$, and $\mu_1,\mu_2,\ldots,\mu_n$ be the zeros of
its matching polynomial. Then,
\begin{equation}
ME(G)=\sum\limits_{i=1}^n|\mu_i|.
\end{equation}
\end{defi}
Moreover, Gutman and Wagner \cite{gutman&Wagner2012} pointed out that the matching
energy is a quantity of relevance for chemical applications. They arrived at the simple
relation: $$TRE(G) = E(G) - ME(G),$$ where TRE($G$) is the so-called ``topological
resonance energy" of $G$. About the chemical applications of matching energy, for more
details see \cite{gutmanmt1975,aihara1976,gutmanmt1977}.

An important tool of graph energy is the Coulson-type integral formula
\cite{gutman1977}\, (with regard to $G$ be a tree $T$):
\begin{equation}\label{cformula}
E(T)=\frac{2}{\pi}\int_0^{\infty}\frac{1}{x^2}\ln\Big[\sum\limits_{k\geq 0}
m(T,k)x^{2k}\Big]dx,
\end{equation}
which is valid for any tree $T$ (or, more generally, for any forest). Being similar
to Eq.(\ref{cformula}), the matching energy also has a beautiful formula as
follows\cite{gutman&Wagner2012}.

\begin{prop}
Let $G$ be a simple graph of order $n$, and $m(G, k)$ be the number of its
$k$-matchings, $k = 0, 1, 2,\ldots,\lfloor\frac{n}{2}\rfloor$. The matching energy
of $G$ is given by
\begin{equation}\label{matchenergyformula}
ME=ME(G)=\frac{2}{\pi}\int^{\infty}_{0}\frac{1}{x^2}\ln\Big[ \sum\limits_{k\geq
0}m(G,k)x^{2k}\Big]dx. \qquad\qquad
\end{equation}
\end{prop}

Combining Eq.(\ref{cformula}) with Eq.(\ref{matchenergyformula}), it immediately
follows that: if $G$ is a forest, then its matching energy coincides with its
energy.

Formula (\ref{cformula}) implies that the energy of a tree is a monotonically
increasing function of any $m(T,k)$. In particular, if $T^{'}$ and $T^{''}$ are two
trees for which $m(T^{'},k)\geq m(T^{''},k)$ holds for all $k\geq1$, then
$E(T^{'})\geq E(T^{''})$. If, in addition, $m(T^{'},k)>m(T^{''},k)$ for at least one
$k$, then $E(T^{'})> E(T^{''})$. Obviously, by Formula (\ref{matchenergyformula})
and the monotonicity of the logarithm function, the result is also valid for $ME$.
Thus, we can define a \emph{quasi-order} ``$\succeq$" as follows: If two graphs
$G_1$ and $G_2$ have the same order and size, then
\begin{equation}\label{quasi-order}
\begin{split}
G_{1}\succeq G_{2}\Longleftrightarrow m(G_1,k)\geq m(G_2,k)\quad \text{for}\, 1\leq k\leq
\left\lfloor\frac{n}{2}\right\rfloor.
\end{split}
\end{equation}
And if $G_{1}\succeq G_{2}$ we say that $G_{1}$ is $m$-\emph{greater than} $G_{2}$ or
$G_{2}$ is $m$-\emph{smaller than} $G_{1}$. If $G_{1}\succeq G_{2}$ and $G_{2}\succeq
G_{1}$, the graphs $G_{1}$ and $G_{2}$ are said to be $m$-\emph{equivalent}, denote it
by $G_{1}\sim G_{2}$. If $G_{1}\succeq G_{2}$, but the graphs $G_{1}$ and $G_{2}$ are
not $m$-equivalent (i.e., there exists some $k$ such that $m(G_{1},k)>m(G_{2},k)$), then
we say that $G_{1}$ is \emph{strictly $m$-greater than} $G_{2}$, write $G_{1}\succ
G_{2}$. If neither $G_{1}\succeq G_{2}$ nor $G_{2}\succeq G_{1}$, the two graphs $G_{1}$
and $G_{2}$ are said to be $m$-\emph{incomparable} and we denote this by $G_{1}\# G_{2}$.

The relation $\sim$ is an equivalence relation in any set of graphs $\gamma$. The corresponding
equivalence classes will be called \emph{matching equivalence classes} (of the set $\gamma$). The
relation $\succeq$ induces a partial ordering of the set $\gamma/\sim$. An equivalence class is
said to be the \emph{greatest} if it is greater than any of other class. A class is \emph{maximal}
if there is no other class greater than it. The graphs belonging to greatest (resp. maximal) matching
equivalence classes will be said to be $m$-\emph{greatest}(resp. $m$-\emph{maximal}) in the set considered.

According to Eq.(\ref{matchenergyformula}) and Eq.(\ref{quasi-order}), we have
$$G_{1}\succeq G_{2} \Longrightarrow ME(G_{1})\geq ME(G_{2})$$ and $$G_{1}\succ G_{2}
\Longrightarrow ME(G_{1})> ME(G_{2}).$$ It follows that the $m$-greatest graphs must
have greatest matching energy, and the $m$-maximal graphs must have greater matching
energy than other graphs not to be $m$-maximal.

%
%We now introduce some elementary notations and terminologies that will be used in the sequel.
%With regard to other notations, the readers are referred to the book \cite{graphBondy2008}.

\begin{figure}[h,t,b,p]
\begin{center}
\includegraphics[scale = 0.7]{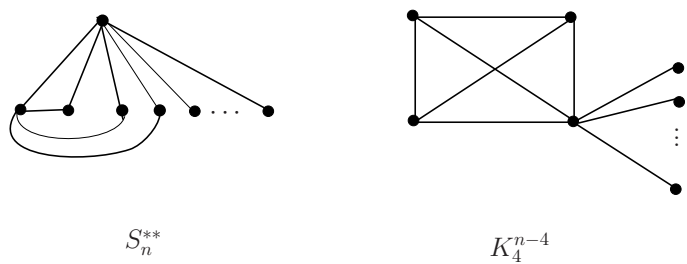}
\caption{Tricyclic graphs with minimal matching energy.}\label{fig1}
\end{center}
\end{figure}

A connected simple graph with $n$ vertices and $n$, $(n+1)$, $(n+2)$ edges are called
\emph{unicyclic, bicyclic, tricyclic graphs}, respectively. Denote by
$\mathscr{B}_n$ the set of all connected bicyclic graphs of order $n$, and by
$\mathscr{T}_{n}$ the set of all connected tricyclic graphs on $n$ vertices.
Let
$S_{n}^{*}$ denote the graph obtained by joining one pendant vertex of $S_n$ to its
other two pendant vertices, respectively. Similarly, let $S_{n}^{**}$ be the graph
obtained by joining one pendent vertex of $S_n$ to its another three pendent
vertices, respectively. Let $K_{4}^{n-4}$ denote the graph obtained by attaching
$(n-4)$ pendent vertices to one of the four vertices of $K_{4}$. Of course,
$S_{n}^{**}, K_{4}^{n-4} \in \mathscr{T}_{n}$ (as shown in Fig. \ref{fig1}).
Denote by $C_{n}$ the cycle graph of order $n$ and $P_{n}$ the path graph of order
$n$, and let $P_n^{k,\ell}$ be the graph obtained by connecting two cycles $C_{k}$ and
$C_{\ell}$ with a path $P_{n-k-\ell}$.

As the research of extremal graph energy is an amusing work (for some newest literatures
see \cite{huo2011a,huo2011,huo2011b,HLS1,HLSW1,LLS1}), the study on extremal matching
energy is also interesting. In \cite{gutman&Wagner2012}, the authors gave some
elementary results on the matching energy and obtained that $ME(S_n^+)\leq ME(G)\leq
ME(C_n)$ for any unicyclic graph $G$, where $S_n^+$ is the graph obtained by adding a
new edge to the star $S_n$. In \cite{ji&li&shi2013}, Ji et al. proved that for $G\in
\mathscr{B}_n$ with $n\geq 10$ and $n=8$,  $ME(S_n^*)\leq ME(G)\leq ME(P_n^{4,n-4})$. In
\cite{jiMa2014}, the authors characterize the connected graphs (and bipartite graph) of
order $n$ having minimum matching energy with $m$ ($n + 2 \leq m \leq 2n-4$) ($n\leq m
\leq 2n-5$) edges. Especially, among all tricyclic graphs of order $n\geq 5$, $ME(G)\geq
ME(S_n^{**})$, with equality if and only if $G\cong S_{n}^{**}$ or $G\cong K_{4}^{n-4}$.
For more results on the matching energy, we refer to \cite{li&yan2014, li&zhou&su2014}.
In this paper, we characterize the graphs with the maximal matching energy among all
tricyclic graphs, and completely determine the tricyclic graphs with the
maximal matching energy.

\section{Main Results}

In the 1980s, Gutman determined the unicyclic \cite{gutman1980}, bicyclic
\cite{gutman1982}, tricyclic \cite{gutmanCvetkovi} graphs with maximal matchings,
i.e., graphs that are extremal with regard to the quasi-ordering $\succeq$. We
introduce the result on tricyclic graphs, which will be used in our proof.

\begin{figure}[h,t,b,p]
\begin{center}
\includegraphics[scale = 0.6]{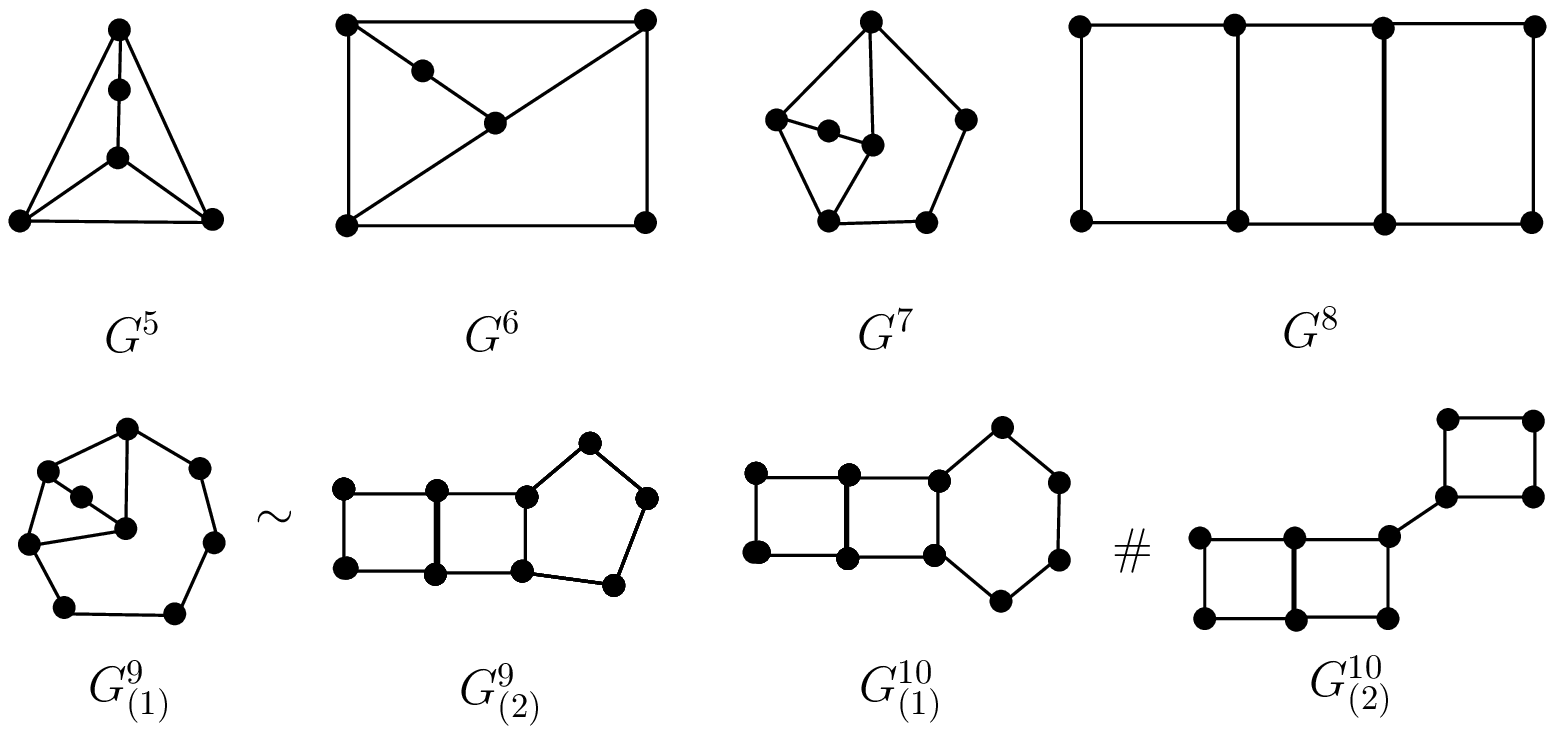}
\includegraphics[scale = 0.6]{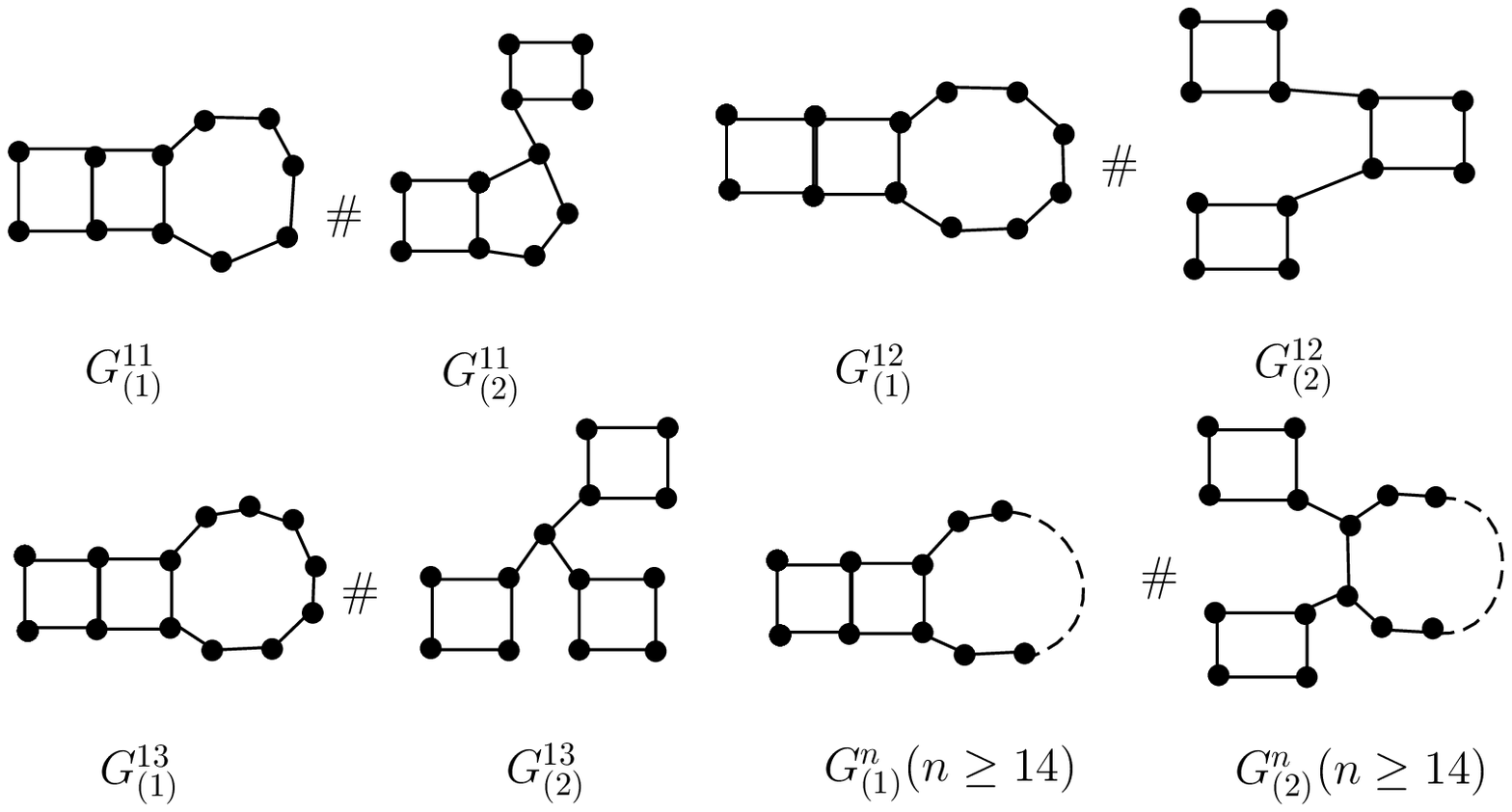}
\caption{The tricyclic graphs with a maximal number of matchings.}\label{fig2}
\end{center}
\end{figure}

\begin{lem}[\cite{gutmanCvetkovi}]\label{maximatching}
In the set of all tricyclic graphs with $n$ vertices ($n\geq 4$) the greatest
matching equivalence class exists only for $n=4, 5, 6, 7, 8 $ and $9$. For $n\geq
10$ there exist two maximal matching equivalence classes. All these equivalence
classes possess a unique element, except for $n=9$, when the number of $m$-greatest
graph is two. The corresponding graphs are presented in Fig. \ref{fig2}.
\end{lem}

Our results are obtained based on the result of Lemma \ref{maximatching}.
\begin{thm}\label{maintri}
Let $G\in \mathscr{T}_n$ with $n\geq 5$. Then for $n=5$, $ME(G)\leq ME(G^{5})$; for
$n=6$, $ME(G)\leq ME(G^{6})$; for $n=7$, $ME(G)\leq ME(G^{7})$; for $n=8$, $ME(G)\leq
ME(G^{8})$; for $n=9$, $ME(G)\leq ME(G_{(1)}^{9})=ME(G_{(2)}^{9})$; for $n=10$,
$ME(G)\leq ME(G_{(2)}^{10})$; for $n=11$, $ME(G)\leq ME(G_{(2)}^{11})$; for $n=12$,
$ME(G)\leq ME(G_{(2)}^{12})$; for $n=13$, $ME(G)\leq ME(G_{(2)}^{13})$; for $n\geq 14$,
$ME(G)\leq ME(G_{(2)}^{n})$, with equality if and only if $G\cong G_{(2)}^{n}$, where
$G^{5}$, $G^{6}$, $G^{7}$, $G^{8}$, $G_{(1)}^{9}$, $G_{(2)}^{9}$, $G_{(2)}^{10}$,
$G_{(2)}^{11}$, $G_{(2)}^{12}$, $G_{(2)}^{13}$, $G_{(2)}^{n}$ are the graphs shown in
Fig. \ref{fig2}.
\end{thm}

We will prove our theorem by using Coulson-type integral formula of matching energy Eq.(\ref{matchenergyformula}),
which is similar as the method to comparing the energies of
two quasi-order incomparable graphs \cite{huo2011a,huo2011,huo2011b,HLS1,HLSW1,LLS1}.
The following lemmas are needed.

%\begin{lem}[\cite{Farrell1979,gutmanmatch1979}]\label{twomindentity}
%Let G be a graph. Then, for any edge e=uv and
%$N(u)=\{v_1(=v),v_2,\ldots,v_t \}$, we have the two identities:
%\begin{equation}\label{removeedge}
%m(G,k)=m(G-uv,k)+m(G-u-v,k-1)
%\end{equation}
%\begin{equation}\label{removevertex}
%m(G,k)=m(G-u,k)+\sum_{i=1}^{t}m(G-u-v_i,k-1).
%\end{equation}
%\end{lem}

%\begin{lem}[\cite{ji&li&shi2013}]\label{bicyclicmini}
%Let $G\in \mathscr{B}_n$ with $n\geq 4$. Then $ME(G)\geq ME(S_n^*)$, with
%equality if and only if $G\cong S_n^*$.
%\end{lem}

%\begin{lem}[\cite{ji&li&shi2013}]\label{bicyclicmaxi}
%Let $G\in \mathscr{B}_n$ with $n\geq 10$ and $n=8$. Then $ME(G)\leq
%ME(P_n^{4,n-4})$, with equality if and only if $G\cong P_n^{4,n-4}$.
%Exceptionally, when $n=9$, $P_n^{4,n-4}$ and $C_n(3,1,n-3)$ are
%matching-equivalent and thus both have maximal ME-values.
%\end{lem}

\begin{lem}[\cite{zorich2002}]\label{realnumber}
For any real number $X>-1$, we have
\begin{equation}\label{realbound}
\frac{X}{1+X}\leq \ln(1+X)\leq X.
\end{equation}
\end{lem}

Let $G$ be a simple graph. Let $e$ be an edge of $G$ connecting the vertices $v_{r}$
and $v_{s}$. By $G(e/j)$ we denote the graph obtained by inserting $j\ (j \geq 0)$
new vertices (of degree two) on the edge $e$. Hence if $G$ has $n$ vertices, then
$G(e/j)$ has $n+j$ vertices; if $j=0$, then $G(e/j)=G$; if $j> 0$, then the vertices
$v_{r}$ and $v_{s}$ are not adjacent in $G(e/j)$. There is such a result on the
number of $k$-matchings of the graph $G(e/j)$.

\begin{lem}[\cite{gutmanCvetkovi}]\label{insertingvertices}
For all $j\geq 0$, $$m(G(e/j+2),k)=m(G(e/j+1),k)+m(G(e/j),k-1).$$
\end{lem}

We will divide Theorem \ref{maintri} into the following two theorems according to the
values of $n$.

\begin{thm}\label{maxi1}
Let $G\in \mathscr{T}_n$ with $n\geq 5$. Then:\\
for $n=5$, $ME(G)\leq ME(G^{5})$; for $n=6$, $ME(G)\leq ME(G^{6})$; for $n=7$,
$ME(G)\leq ME(G^{7})$; for $n=8$, $ME(G)\leq ME(G^{8})$; for $n=9$, $ME(G)\leq
ME(G_{(1)}^{9})=ME(G_{(2)}^{9})$; for $n=10$, $ME(G)\leq ME(G_{(2)}^{10})$; for
$n=11$, $ME(G)\leq ME(G_{(2)}^{11})$; for $n=12$, $ME(G)\leq ME(G_{(2)}^{12})$; for
$n=13$, $ME(G)\leq ME(G_{(2)}^{13})$, where $G^{5}$, $G^{6}$, $G^{7}$, $G^{8}$,
$G_{(1)}^{9}$, $G_{(2)}^{9}$, $G_{(2)}^{10}$, $G_{(2)}^{11}$, $G_{(2)}^{12}$,
$G_{(2)}^{13}$ are the graphs shown in Fig. \ref{fig2}. In each case, the equality
holds if and only if $G$ is isomorphic to the corresponding graph with maximal
matching energy.
\end{thm}

\pf Let $G$ be a graph in $\mathscr{T}_n$ with $n$ vertices.

For $n=5, 6, 7, 8$, by Lemma \ref{maximatching}, $G^{n}$ is the $m$-greatest graph. We have known that the $m$-greatest graphs
must have greatest matching energy. Hence if $G\ncong G^{n}$, then $ME(G)< ME(G^{n})$.

When $n=9$, $G_{(1)}^{9}$ and $G_{(2)}^{9}$ are $m$-equivalent, that is,
$m(G_{(1)}^{9},k)=m(G_{(2)}^{9},k)$ for all $k$. Then by Eq.(\ref{matchenergyformula}),
we have $ME(G_{(1)}^{9})=ME(G_{(2)}^{9})$. Moreover, if
$G\ncong G_{(1)}^{9}$ and $G\ncong G_{(2)}^{9}$, then $ME(G)<
ME(G_{(1)}^{9})=ME(G_{(2)}^{9})$ since $(G_{(1)}^{9}\sim G_{(2)}^{9})\succ G$ by
Lemma \ref{maximatching}.

When $n=10$, both $G_{(1)}^{10}$ and $G_{(2)}^{10}$ are $m$-maximal. Thus, if
$G\ncong G_{(1)}^{10}$ and $G\ncong G_{(2)}^{10}$, then $ME(G)< ME(G_{(1)}^{10})$
as well as $ME(G)<ME(G_{(2)}^{10})$. In addition, we have $m(G^{10}_{(1)},0)=1$,
$m(G^{10}_{(1)},1)=12$, $m(G^{10}_{(1)},2)=48$, $m(G^{10}_{(1)},3)=76$,
$m(G^{10}_{(1)},4)=42$, $m(G^{10}_{(1)},5)=5$ and $m(G^{10}_{(2)},0)=1$,
$m(G^{10}_{(2)},1) =12$, $m(G^{10}_{(2)},2)=48$, $m(G^{10}_{(2)},3)=75$,
$m(G^{10}_{(2)},4)=42$, $m(G^{10}_{(2)},5)=6$. Make use of Eq.(\ref{matchenergyformula}),
by computer-aided calculations, we get
$ME(G^{10}_{(1)})=13.8644$ and $ME(G^{10}_{(2)})=13.9042$. Therefore, $ME(G)<
ME(G_{(1)}^{10})< ME(G_{(2)}^{10})$.

For $n=11, 12, 13$, both $G^{n}_{(1)}$ and $G^{n}_{(2)}$ are $m$-maximal. Similarly,
by the help of computer, we get $ME(G^{11}_{(1)})=14.9384$,
$ME(G^{11}_{(2)})=14.9466$, $ME(G^{12}_{(1)})=16.3946$, $ME(G^{12}_{(2)})=16.5052$,
$ME(G^{13}_{(1)})=17.5097$, $ME(G^{13}_{(2)})=17.5678$, respectively. Therefore, if
$G\ncong G_{(2)}^{n}$, then we have $ME(G)\leq ME(G_{(1)}^{n})< ME(G_{(2)}^{n})$.

The proof of the theorem is complete.  \qed

\begin{thm}\label{maxi2}
Let $G\in \mathscr{T}_n$ with $n\geq 14$. Then $ME(G)\leq ME(G_{(2)}^{n})$, with
equality if and only if $G\cong G_{(2)}^{n}$, where $G_{(2)}^{n}$ is the graph shown
in Fig. \ref{fig2}.
\end{thm}

\pf By Lemma \ref{maximatching}, both $G^{n}_{(1)}$ and $G^{n}_{(2)}$ are
$m$-maximal. The $m$-maximal graphs must have greater matching energy than other
graphs not to be $m$-maximal. Thus, if $G\ncong G_{(1)}^{n}$ and $G\ncong
G_{(2)}^{n}$, then $ME(G)< ME(G_{(1)}^{n})$ and $ME(G)< ME(G_{(2)}^{n})$. It is
sufficient to show that $ME(G_{(1)} ^{n})< ME(G_{(2)}^{n})$. We will make full use
of the definition of matching polynomial and Eq.(\ref{matchenergyformula}).

Assume that $|G(e/j+2)|=n$, then$|G(e/j+1)|=n-1$ and $|G(e/j)|=n-2$. According to
Lemma \ref{insertingvertices}, we have
\begin{equation*}
\begin{split}
\alpha(G(e/j+2),x)&=\sum\limits_{k\geq 0}(-1)^km(G(e/j+2),k)x^{n-2k}\\
                  &=\sum\limits_{k\geq 0}(-1)^km(G(e/j+1),k)x^{n-2k}+\sum\limits_{k\geq 0}(-1)^km(G(e/j),k-1)x^{n-2k}\\
                  &=x\sum\limits_{k\geq 0}(-1)^km(G(e/j+1),k)x^{(n-1)-2k}\\
                  &-\sum\limits_{k\geq 0}(-1)^{k-1}m(G(e/j),k-1)x^{(n-2)-2(k-1)}\\
                  &=x\alpha(G(e/j+1),x)-\alpha(G(e/j),x).\\
\end{split}
\end{equation*}

\begin{figure}[h,t,b,p]
\begin{center}
\includegraphics[scale = 0.7]{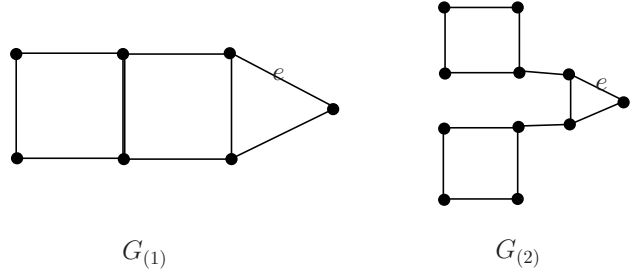}
\caption{The fundamental graphs for constructing $G_{(1)}^{n}$ and
$G_{(2)}^{n}$.}\label{fig7}
\end{center}
\end{figure}

By the definition of $G(e/j)$, clearly, $G_{(1)}^{n}=G_{(1)}(e/n-7)$ and
$G_{(2)}^{n}=G_{(2)}(e/n-11)$, where $G_{(1)}$ and $G_{(2)}$ are the graphs shown in
Fig. \ref{fig7}.
Therefore, both $\alpha(G_{(1)}^{n},x)$ and $\alpha(G_{(2)}^{n},x)$ satisfy the
recursive formula $$f(n,x)=xf(n-1,x)-f(n-2,x).$$ The general solution of this linear
homogeneous recurrence relation is $$f(n,x)=
C_{1}(x)(Y_{1}(x))^{n}+C_{2}(x)(Y_{2}(x))^{n},$$ where
$Y_{1}(x)=\frac{x+\sqrt{x^{2}-4}}{2}$, $Y_{2}(x)=\frac{x-\sqrt{x^{2}-4}}{2}$. By
some elementary calculations, we can easily obtain the values of $C_{i}(x)\ (i=1,
2)$ as follows.

In the following, we first consider $\alpha(G_{(1)}^{n},x)$. It is easy to calculate
the number of $k$-matchings of $G_{(1)}$ and $G_{(1)}(e/1)$: $m(G_{(1)},0)=1$,
$m(G_{(1)},1)=9$, $m(G_{(1)},2)=21$, $m(G_{(1)},3)=11$, $m(G_{(1)},k)=0$ for $k\geq
4$; $m(G_{(1)}(e/1),0)=1$, $m(G_{(1)}(e/1),1)=10$, $m(G_{(1)}(e/1),2)=29$,
$m(G_{(1)}(e/1),3) =26$, $m(G_{(1)}(e/1),4)=5$, $m(G_{(1)}(e/1),k)=0$ for $k\geq 5$.
Then by Lemma \ref{insertingvertices}, we can calculate the values of
$m(G_{(1)}(e/j),k)$ for all $j\geq 2$ and $k\geq 0$. Thus, take the initial values
as:
\begin{equation*}
\begin{split}
\alpha(G_{(1)}(e/4),x)&=x^{11}-13x^{9}+59x^{7}-114x^{5}+89x^{3}-21x\\
                &=C_{1}(x)(Y_{1}(x))^{11}+C_{2}(x)(Y_{2}(x))^{11};\\
\alpha(G_{(1)}(e/5),x)&=x^{12}-14x^{10}+71x^{8}-162x^{6}+165x^{4}-63x^{2}+5\\
                &=C_{1}(x)(Y_{1}(x))^{12}+C_{2}(x)(Y_{2}(x))^{12}.\\
\end{split}
\end{equation*}
It is easy to check that $Y_{1}(x)+Y_{2}(x)=x$ and $Y_{1}(x)\cdot Y_{2}(x)=1$.
Therefore, by solving the two equalities above, we get
$$C_{1}(x)=\frac{Y_{1}(x)\alpha(G_{(1)}(e/5),x)-\alpha(G_{(1)}(e/4),x)}{(Y_{1}(x))^{13}-(Y_{1}(x))^{11}}$$
and
$$C_{2}(x)=\frac{Y_{2}(x)\alpha(G_{(1)}(e/5),x)-\alpha(G_{(1)}(e/4),x)}{(Y_{2}(x))^{13}-(Y_{2}(x))^{11}}.$$
Define
$$A_{1}(x)=\frac{Y_{1}(x)\alpha(G_{(1)}(e/5),x)-\alpha(G_{(1)}(e/4),x)}{(Y_{1}(x))^{13}-(Y_{1}(x))^{11}}.$$
$$A_{2}(x)=\frac{Y_{2}(x)\alpha(G_{(1)}(e/5),x)-\alpha(G_{(1)}(e/4),x)}{(Y_{2}(x))^{13}-(Y_{2}(x))^{11}}.$$
Then we have $\alpha(G_{(1)}^{n},x)=A_{1}(x)(Y_{1}(x))^{n}+A_{2}(x)(Y_{2}(x))^{n}$.

Now we consider $\alpha(G_{(2)}^{n},x)$. Similarly, we get: $m(G_{(2)},0)=1$,
$m(G_{(2)},1)=13$, $m(G_{(2)},2)=59$, $m(G_{(2)},3)=112$, $m(G_{(2)},4)=84$,
$m(G_{(2)},5)=20$, $m(G_{(2)},k)=0$ for $k\geq 6$; $m(G_{(2)}(e/1),0)=1$,
$m(G_{(2)}(e/1),1)=14$, $m(G_{(2)}(e/1),2)=71$, $m(G_{(2)}(e/1),3)=161$,
$m(G_{(2)}(e/1),4)=164$, $m(G_{(2)}(e/1),5)=68$, $m(G_{(2)}(e/1),6)=8$,
$m(G_{(2)}(e/1),k)=0$ for $k\geq 7$. Then calculate the values of
$m(G_{(2)}(e/j),k)$ for all $j\geq 2$ and $k\geq 0$ by using Lemma
\ref{insertingvertices}. We can then take the initial values as:
\begin{equation*}
\begin{split}
\alpha(G_{(2)},x)&=x^{11}-13x^{9}+59x^{7}-112x^{5}+84x^{3}-20x\\
                 &=C_{1}(x)(Y_{1}(x))^{11}+C_{2}(x)(Y_{2}(x))^{11};\\
\alpha(G_{(2)}(e/1),x)&=x^{12}-14x^{10}+71x^{8}-161x^{6}+164x^{4}-68x^{2}+8\\
                      &=C_{1}(x)(Y_{1}(x))^{12}+C_{2}(x)(Y_{2}(x))^{12}.\\
\end{split}
\end{equation*}
Therefore, we obtain that:
$$C_{1}(x)=\frac{Y_{1}(x)\alpha(G_{(2)}(e/1),x)-\alpha(G_{(2)},x)}{(Y_{1}(x))^{13}-(Y_{1}(x))^{11}}$$
and
$$C_{2}(x)=\frac{Y_{2}(x)\alpha(G_{(2)}(e/1),x)-\alpha(G_{(2)},x)}{(Y_{2}(x))^{13}-(Y_{2}(x))^{11}}.$$
Define
$$B_{1}(x)=\frac{Y_{1}(x)\alpha(G_{(2)}(e/1),x)-\alpha(G_{(2)},x)}{(Y_{1}(x))^{13}-(Y_{1}(x))^{11}}.$$
$$B_{2}(x)=\frac{Y_{2}(x)\alpha(G_{(2)}(e/1),x)-\alpha(G_{(2)},x)}{(Y_{2}(x))^{13}-(Y_{2}(x))^{11}}.$$
Then we have $\alpha(G_{(2)}^{n},x)=B_{1}(x)(Y_{1}(x))^{n}+B_{2}(x)(Y_{2}(x))^{n}.$

From the expression of $\alpha(G,x)$, we have
\begin{equation*}
\begin{split}
\alpha(G,ix)&=\sum\limits_{k\geq 0}(-1)^km(G,k)(ix)^{n-2k}=i^{n}\sum\limits_{k\geq
0}m(G,k)x^{n-2k}
            =(ix)^{n}\sum\limits_{k\geq 0}m(G,k)x^{-2k},
\end{split}
\end{equation*}
where $i^{2}=-1$. Thus, by Eq.(\ref{matchenergyformula}), we get
\begin{equation}\label{MEdifference}
\begin{split}
ME(G_{(1)}^{n})-ME(G_{(2)}^{n})&=\frac{2}{\pi}\int^{\infty}_{0}\frac{1}{x^2}\ln\Big[\sum\limits_{k\geq 0}m(G_{(1)}^{n},k)x^{2k}\Big]dx\\
                               &-\frac{2}{\pi}\int^{\infty}_{0}\frac{1}{x^2}\ln\Big[\sum\limits_{k\geq 0}m(G_{(2)}^{n},k)x^{2k}\Big]dx\\
                               &=\frac{2}{\pi}\int^{\infty}_{0}\frac{1}{x^2}\ln\frac{\sum\limits_{k\geq 0}m(G_{(1)}^{n},k)x^{2k}}{\sum\limits_{k\geq 0}m(G_{(2)}^{n},k)x^{2k}}dx\\
                               &=\frac{2}{\pi}\int^{\infty}_{0}\ln\frac{\sum\limits_{k\geq 0}m(G_{(1)}^{n},k)x^{-2k}}{\sum\limits_{k\geq 0}m(G_{(2)}^{n},k)x^{-2k}}dx\\
                               &=\frac{2}{\pi}\int^{\infty}_{0}\ln\frac{\alpha(G_{(1)}^{n},ix)}{\alpha(G_{(2)}^{n},ix)}dx\\
                               &=\frac{2}{\pi}\int^{\infty}_{0}\ln\frac{A_{1}(ix)(Y_{1}(ix))^{n}+A_{2}(ix)(Y_{2}(ix))^{n}}
                               {B_{1}(ix)(Y_{1}(ix))^{n}+B_{2}(ix)(Y_{2}(ix))^{n}}dx.
\end{split}
\end{equation}

By the definition of $Y_{1}(x)$ and $Y_{2}(x)$, we have
$Y_{1}(ix)=\frac{x+\sqrt{x^{2}+4}}{2}i$ and $Y_{2}(ix)=\frac{x-\sqrt{x^{2}+4}}{2}i$.
Now we define $Z_{1}(x)=-iY_{1}(x)=\frac{x+\sqrt{x^{2}+4}}{2}$, $Z_{2}(x)=
-iY_{2}(x)=\frac{x-\sqrt{x^{2}+4}}{2}$, and
$$f_{1}=i\alpha(G_{(1)}(e/4),ix)=x^{11}+13x^{9}+59x^{7}+114x^{5}+89x^{3}+21x$$
$$f_{2}=\alpha(G_{(1)}(e/5),ix)=x^{12}+14x^{10}+71x^{8}+162x^{6}+165x^{4}+63x^{2}+5$$
$$g_{1}=i\alpha(G_{(2)},ix)=x^{11}+13x^{9}+59x^{7}+112x^{5}+84x^{3}+20x$$
$$g_{2}=\alpha(G_{(2)}(e/1),ix)=x^{12}+14x^{10}+71x^{8}+161x^{6}+164x^{4}+68x^{2}+8.$$
Then we have $Y_{1}(ix)=iZ_{1}(x)$ and $Y_{2}(ix)=iZ_{2}(x)$. Moreover, It follows that
$$A_{1}(ix)=\frac{iZ_{1}(x)f_{2}+if_{1}}{(iZ_{1}(x))^{13}-(iZ_{1}(x))^{11}}=\frac{Z_{1}(x)f_{2}+f_{1}}{(Z_{1}(x))^{11}((Z_{1}(x))^{2}+1)}$$
$$A_{2}(ix)=\frac{iZ_{2}(x)f_{2}+if_{1}}{(iZ_{2}(x))^{13}-(iZ_{2}(x))^{11}}=\frac{Z_{2}(x)f_{2}+f_{1}}{(Z_{2}(x))^{11}((Z_{2}(x))^{2}+1)}$$
$$B_{1}(ix)=\frac{iZ_{1}(x)g_{2}+ig_{1}}{(iZ_{1}(x))^{13}-(iZ_{1}(x))^{11}}=\frac{Z_{1}(x)g_{2}+g_{1}}{(Z_{1}(x))^{11}((Z_{1}(x))^{2}+1)}$$
$$B_{2}(ix)=\frac{iZ_{2}(x)g_{2}+ig_{1}}{(iZ_{2}(x))^{13}-(iZ_{2}(x))^{11}}=\frac{Z_{2}(x)g_{2}+g_{1}}{(Z_{2}(x))^{11}((Z_{2}(x))^{2}+1)}.$$
Note that $Y_{1}(ix)\cdot Y_{2}(ix)=1$, $Z_{1}(x)\cdot Z_{2}(x)=-1$, $Z_{1}(x)+Z_{2}(x)=x$ and $Z_{1(x)}-Z_{2}(x)=\sqrt{x^{2}+4}$.
We will distinguish with two cases.

\noindent{\bf Case 1.} $n$ is odd.

Now we have
\begin{equation*}
\begin{split}
&\ln\frac{A_{1}(ix)(Y_{1}(ix))^{n+2}+A_{2}(ix)(Y_{2}(ix))^{n+2}}{B_{1}(ix)(Y_{1}(ix))^{n+2}+B_{2}(ix)(Y_{2}(ix))^{n+2}}
-\ln\frac{A_{1}(ix)(Y_{1}(ix))^{n}+A_{2}(ix)(Y_{2}(ix))^{n}}{B_{1}(ix)(Y_{1}(ix))^{n}+B_{2}(ix)(Y_{2}(ix))^{n}}\\[2mm]
&=\ln\left(1+\frac{K_{0}(x)}{H_{0}(n,x)}\right),
\end{split}
\end{equation*}
where
\begin{equation*}
\begin{split}
K_{0}(x)&=(A_{1}(ix)B_{2}(ix)-A_{2}(ix)B_{1}(ix))((Y_{1}(ix))^{2}-(Y_{2}(ix))^{2})
        =(f_{2}g_{1}-f_{1}g_{2})x\\
        &=-x^{18}-19x^{16}-146x^{14}-588x^{12}-1342x^{10}-1750x^{8}-1253x^{6}-460x^{4}-68x^{2},
\end{split}
\end{equation*}
and
\begin{equation*}
\begin{split}
H_{0}(n,x)&=(A_{1}(ix)(Y_{1}(ix))^{n}+A_{2}(ix)(Y_{2}(ix))^{n})(B_{1}(ix)(Y_{1}(ix))^{n+2}+B_{2}(ix)(Y_{2}(ix))^{n+2})\\
          &=\alpha(G_{(1)}^{n},ix)\cdot \alpha(G_{(2)}^{n+2},ix)\\
          &=\left(i^{n}\sum\limits_{k\geq 0}m(G_{(1)}^{n},k)x^{n-2k}\right)\left(i^{n+2}\sum\limits_{k\geq 0}m(G_{(2)}^{n+2},k)x^{n+2-2k}\right)\\
          &=i^{2n+2}\left(\sum\limits_{k\geq 0}m(G_{(1)}^{n},k)x^{n-2k}\right)\left(\sum\limits_{k\geq 0}m(G_{(2)}^{n+2},k)x^{n+2-2k}\right).
\end{split}
\end{equation*}

Obviously, $K_{0}(x)<0$. Moreover, since $n$ is odd, we have $i^{2n+2}=1$, it
follows that $H_{0}(n,x)$ is a polynomial such that each term is of positive even
degree of $x$ and all coefficients are positive, i.e., $H_{0}(n,x)>0$. Hence,
$\frac{K_{0}(x)} {H_{0}(n,x)}<0$, which deduces that
$\ln(1+\frac{K_{0}(x)}{H_{0}(n,x)})<\ln1=0$ for $x>0$ and odd $n$. So, the integrand
of Eq.(\ref{MEdifference}) is monotonically decreasing on $n$. Therefore, for
$n\geq 14$,
$$\int^{\infty}_{0}\ln\frac{\alpha(G_{(1)}^{n},ix)}{\alpha(G_{(2)}^{n},ix)}dx\leq \int^{\infty}_{0}\ln\frac{\alpha(G_{(1)}^{15},ix)}
{\alpha(G_{(2)}^{15},ix)}dx=\int^{\infty}_{0}\ln\frac{\alpha(G_{(1)}(e/8),ix)}{\alpha(G_{(2)}(e/4),ix)}dx.$$
By computer-aided calculations, we get $ME(G_{(1)}(e/8))=20.0728$ and $ME(G_{(2)}(e/4))=20.1086$. And then
$$\int^{\infty}_{0}\ln\frac{\alpha(G_{(1)}(e/8),ix)}{\alpha(G_{(2)}(e/4),ix)}dx=\frac{\pi}{2}\left(ME(G_{(1)}(e/8))-ME(G_{(2)}(e/4))\right)=-0.05639<0.$$
So $\int^{\infty}_{0}\ln\frac{\alpha(G_{(1)}^{n},ix)}{\alpha(G_{(2)}^{n},ix)}dx<0$.
That is, $$ME(G_{(1)}^{n})-ME(G_{(2)}^{n})=
\frac{2}{\pi}\int^{\infty}_{0}\ln\frac{\alpha(G_{(1)}^{n},ix)}{\alpha(G_{(2)}^{n},ix)}dx<0.$$
Therefore, $ME(G_{(1)}^{n})< ME(G_{(2)}^{n})$ when $n$ is odd.

\noindent{\bf Case 2.} $n$ is even.

Since $x>0$, when $n\longrightarrow \infty$, we have
$$\frac{A_{1}(ix)(Y_{1}(ix))^{n}+A_{2}(ix)(Y_{2}(ix))^{n}}{B_{1}(ix)(Y_{1}(ix))^{n}+B_{2}(ix)(Y_{2}(ix))^{n}}
\longrightarrow \frac{A_{1}(ix)}{B_{1}(ix)}.$$
Then we have
$$\ln\frac{A_{1}(ix)(Y_{1}(ix))^{n}+A_{2}(ix)(Y_{2}(ix))^{n}}{B_{1}(ix)(Y_{1}(ix))^{n}+B_{2}(ix)(Y_{2}(ix))^{n}}
-\ln\frac{A_{1}(ix)}{B_{1}(ix)}\\
=\ln\left(1+\frac{K_{1}(n,x)}{H_{1}(n,x)}\right),$$ where
\begin{equation*}
\begin{split}
K_{1}(n,x)&=A_{2}(ix)\cdot B_{1}(ix)\cdot (Y_{2}(ix))^{n}-A_{1}(ix)\cdot B_{2}(ix)\cdot (Y_{2}(ix))^{n}\\
          &=\frac{(f_{2}g_{1}-f_{1}g_{2})(Z_{2}(x))^{n}}{\sqrt{x^{2}+4}}\cdot i^{n}\\
          &=\frac{(-x^{17}-19x^{15}-146x^{13}-588x^{11}-1342x^{9}-1750x^{7}-1253x^{5}-460x^{3}-68x)}{\sqrt{x^{2}+4}}\\
          &~~~~~\cdot (Z_{2}(x))^{n}\cdot i^{n},
\end{split}
\end{equation*}
and
\begin{equation*}
\begin{split}
H_{1}(n,x)&=A_{1}(ix)\cdot B_{1}(ix)\cdot (Y_{1}(ix))^{n}+A_{1}(ix)\cdot B_{2}(ix)\cdot (Y_{2}(ix))^{n}\\
          &=A_{1}(ix)(B_{1}(ix)\cdot (Y_{1}(ix))^{n}+B_{2}(ix)\cdot (Y_{2}(ix))^{n})
          =A_{1}(ix)\alpha(G_{(2)}^{n},ix)\\
          &=i^{n}\cdot \frac{Z_{1}(x)f_{2}+f_{1}}{(Z_{1}(x))^{11}((Z_{1}(x))^{2}+1)}\cdot \sum\limits_{k\geq 0}m(G_{(2)}^{n},k)x^{n-2k}.
\end{split}
\end{equation*}

Since $n$ is even, $(Z_{2}(x))^{n}>0$. Hence $K_{1}(n,x)/i^{n}$ is a polynomial of
$x$ with all coefficients being negative, namely, we always have $K_{1}(n,x)/i^{n}<
0$. On the other hand, since $x>0$, we have $Z_{1}(x)=\frac{x+\sqrt
{x^{2}+4}}{2}>0$, $f_{1}>0$, $f_{2}>0$ and $m(G_{(2)}^{n},k)> 0$ for all $0\leq k\leq \lfloor\frac{n}{2}\rfloor$.
Hence $H_{1}(n,x)/i^{n}$ is a polynomial of $x$ such that all the coefficients are
positive. Therefore, $\frac{K_{1}(n,x)}{H_{1}(n,x)}< 0$ for all $x>0$ and even $n$.
Then $\ln(1+\frac{K_{1}(n,x)}{H_{1}(n,x)})<\ln1=0$, i.e.,
$$\ln\frac{A_{1}(ix)(Y_{1}(ix))^{n}+A_{2}(ix)(Y_{2}(ix))^{n}}{B_{1}(ix)(Y_{1}(ix))^{n}+B_{2}(ix)(Y_{2}(ix))^{n}}<
\ln\frac{A_{1}(ix)}{B_{1}(ix)}.$$ Thus, we have proved that the integrand of Eq.(\ref{MEdifference})
is less than the corresponding limit function when $n$ is even.
Furthermore, since
$$1+\frac{A_{1}(ix)-B_{1}(ix)}{B_{1}(ix)}=\frac{A_{1}(ix)}{B_{1}(ix)}=\frac{Z_{1}(x)f_{2}+f_{1}}
{Z_{1}(x)g_{2}+g_{1}}>0,$$ we have $\frac{A_{1}(ix)-B_{1}(ix)}{B_{1}(ix)}>-1$. Then
by Lemma \ref{realnumber}, $\ln\frac{A_{1} (ix)}{B_{1}(ix)}\leq
\frac{A_{1}(ix)-B_{1}(ix)}{B_{1}(ix)}$. By some computer-aided calculations, we
obtain that $\int^{\infty}_{0}\frac{A_{1}(ix)-B_{1}(ix)}{B_{1}(ix)}dx=-0.09693$. It
means that
$$\int^{\infty}_{0}\ln\frac{A_{1}(ix)}{B_{1}(ix)}dx\leq \int^{\infty}_{0}\frac{A_{1}(ix)-B_{1}(ix)}{B_{1}(ix)}dx<0.$$
Thus,
\begin{equation*}
\begin{split}
\frac{\pi}{2}(ME(G_{(1)}^{n})-ME(G_{(2)}^{n}))&=\int^{\infty}_{0}\ln\frac{A_{1}(ix)(Y_{1}(ix))^{n}+A_{2}(ix)(Y_{2}(ix))^{n}}{B_{1}(ix)(Y_{1}(ix))^{n}+B_{2}(ix)(Y_{2}(ix))^{n}}dx\\[2mm]
&<\int^{\infty}_{0}\ln\frac{A_{1}(ix)}{B_{1}(ix)}dx<0,
\end{split}
\end{equation*}
i.e., $ME(G_{(1)}^{n})<ME(G_{(2)}^{n})$ when $n$ is even.

Therefore, for all $n\geq 14$, we can always show that
$$ME(G_{(1)}^{n})<ME(G_{(2)}^{n}),$$ the proof is thus complete. \qed

{\bf Acknowledgement.} The authors are very grateful to Professor Ivan Gutman for
providing us with reference \cite{gutmanCvetkovi}. The authors are supported by NSFC,
PCSIRT, China Postdoctoral Science Foundation (2014M551015) and China Scholarship
Council.

%\newpage
%\begin{appendices}
%
%
%
%\begin{figure}[h,t,b,p]
%\begin{center}
%\includegraphics[scale = 0.7]{7cycles.eps}
%\caption{One possible case for the arrangement of seven cycles in $G$.}\label{fig3}
%\end{center}
%\end{figure}
%
%\begin{figure}[h,t,b,p]
%\begin{center}
%\includegraphics[scale = 0.7]{6cycles.eps}
%\caption{Three possible cases for the arrangement of six cycles in $G$.}\label{fig4}
%\end{center}
%\end{figure}
%
%\begin{figure}[h,t,b,p]
%\begin{center}
%\includegraphics[scale = 0.7]{4cycles.eps}
%\caption{Four possible cases for the arrangement of four cycles in $G$.}\label{fig5}
%\end{center}
%\end{figure}
%
%\begin{figure}[h,t,b,p]
%\begin{center}
%\includegraphics[scale = 0.7]{3cycles.eps}
%\caption{Seven possible cases for the arrangement of three cycles in $G$.}\label{fig6}
%\end{center}
%\end{figure}
%
%\begin{figure}[h,t,b,p]
%\begin{center}
%\includegraphics[scale = 0.7]{somegraphs1.eps}
%\includegraphics[scale = 0.7]{somegraphs2.eps}
%\caption{Some graphs needed in the proof of Theorem \ref{mini}.}\label{fig8}
%\end{center}
%\end{figure}
%
%\begin{figure}[h,t,b,p]
%\begin{center}
%\includegraphics[scale = 0.7]{somegraphs3.eps}
%\includegraphics[scale = 0.7]{somegraphs4.eps}
%\caption{Some graphs needed in the proof of Theorem \ref{mini}.}\label{fig9}
%\end{center}
%\end{figure}
%
%\end{appendices}

\end{document}